\documentclass[11pt]{article}
\usepackage{amsfonts}
\usepackage{mathrsfs}
\usepackage{amsthm}
 \usepackage{graphicx} 
\usepackage{amsmath} 
\usepackage{amssymb} 
\usepackage{enumerate}

\usepackage{epsfig}

\newtheorem{qe}{\hskip\parindent Question}
\newtheorem{thm}{\hskip\parindent Theorem}
\newtheorem{lem}{\hskip\parindent Lemma}
\newtheorem{prop}{\hskip\parindent Proposition}

\newtheorem{rem}{\hskip\parindent Remark}

\numberwithin{equation}{section}

\def\Re {\mathop{\rm Re}\nolimits}
\def\Im {\mathop{\rm Im}\nolimits}

\begin{document}
\title{  Asymptotics of Landau constants with optimal error bounds}
\author{Yutian Li$^a$, Saiyu Liu$^b$, Shuaixia Xu$^c$ and Yuqiu Zhao$^d$\footnote{Corresponding author.
 {\it{E-mail
address:}} {stszyq@mail.sysu.edu.cn} }}
 \date{ {\it{$^a$Institute of Computational and Theoretical Studies, and Department of Mathematics,   Hong Kong Baptist University, Kowloon, Hong Kong}}\\
 {\it{$^b$School of  Mathematics and Computational Science, Hunan  University of Science and Technology, Xiangtan 411201, Hunan, China}} \\
 {\it{$^c$Institut Franco-Chinois de    l'Energie Nucl\'{e}aire,   Sun Yat-sen University, GuangZhou
510275,  China}}\\
 {\it{$^d$Department of Mathematics, Sun Yat-sen University, GuangZhou
510275, China}}
}

\maketitle

\vskip .3cm
 \noindent {\bf{Abstract: }}
We study the  asymptotic expansion for the Landau constants $G_n$
$$\pi G_n\sim   \ln  N  + \gamma+4\ln 2  + \sum_{s=1}^\infty \frac { \beta_{2s}}{ N^{2s}},~~n\rightarrow \infty, $$
where $N=n+3/4$, $\gamma=0.5772\cdots$ is Euler's constant, and   $(-1)^{s+1}\beta_{2s}$ are positive rational numbers, given  explicitly in an iterative manner. We show that the error due to truncation is  bounded in absolute value  by, and of the same sign as,  the first neglected term for all nonnegative $n$. Consequently, we obtain optimal  sharp bounds up to arbitrary orders of the form
  $$  \ln N+\gamma+4\ln 2+\sum_{s=1}^{2m}\frac{\beta_{2s}}{N^{2s}}<     \pi G_n <  \ln N+\gamma+4\ln 2+\sum_{s=1}^{2k-1}\frac{\beta_{2s}}{N^{2s}}$$
 for all $n=0,1,2,\cdots$, $m=1,2,\cdots$, and $k=1,2,\cdots$.

 The results are proved by approximating  the coefficients $\beta_{2s}$  with  the Gauss hypergeometric functions involved, and by
 using the second order difference equation satisfied by $G_n$, as well as an integral representation of the constants $\rho_k=(-1)^{k+1}\beta_{2k}/(2k-1)!$.
 \vskip .3cm

\noindent {\it{MSC2010:}} 39A60; 41A60; 41A17; 33C05

\vskip .3cm \noindent {\it {Keywords: }} Landau constants; second-order linear difference equation;  asymptotic expansion; sharper bound. \vskip.3cm


 \section{Introduction and statement of results} \indent\setcounter{section} {1}
\setcounter{equation} {0} \label{sec:1}

A  century ago, it was
shown by Landau \cite{Landau} that if a function
$f(z)$ is analytic in the unit disc, such that $|f(z)|<1$,   with
 the Maclaurin expansion
$$f(z)=a_0+a_1 z+a_2 z^2+\cdots+ a_n z^n+\cdots ,~~~|z|<1,$$
then   it holds
$$|a_0+a_1+a_2+\cdots +a_n|\leq G_n,~~~n=0, 1,2,\cdots,$$
where $G_0=1$ and
\begin{equation}\label{(1.1)}
G_n=1+ \left(\frac 1 2\right)^2 +\left ( \frac {1\cdot 3}{2\cdot 4} \right )^2 +\cdots +\left ( \frac{1\cdot 3\cdots (2n-1)}{2\cdot 4\cdots (2n)}\right )^2
\end{equation}
for $n=1,2,\cdots$,
and the equal sign can be attained for each $n$.  The constants   $G_n$ are termed Landau's constants; see, e.g., Watson \cite{Watson}.

Efforts have been made to approximate these constants from the very beginning. Indeed,
 Landau himself \cite{Landau} has worked out the large-$n$ behavior
\begin{equation*}\label{(1.2)}
G_n\sim \frac 1 \pi \ln n,~~~\mbox{as}~~n\rightarrow \infty; \end{equation*}
see also Watson \cite{Watson}.

Since then, the approximation of $G_n$ goes to two related directions. One is to find sharper bounds of $G_n$ for all positive integers $n$, and the other is
to obtain  large-$n$ asymptotic approximations for the constants.

\subsection{Sharper bounds}

Many authors have worked on the sharp bounds   of $G_n$.
 For example, in 1982, Brutman \cite{Brutman}  obtains
\begin{equation*}\label{(1.3)} 1+  \pi^{-1} \ln(n+1) \leq G_n< 1.0663+    \pi^{-1} \ln(n+1),~~n=0,1,2,\cdots. \end{equation*}
The result is improved in 1991 by Falaleev \cite{Falaleev} to give
\begin{equation*}\label{(1.4)} 1.0662+\pi^{-1} \ln(n+0.75)< G_n \leq 1.0916 +\pi^{-1} \ln(n+0.75),~~n=0,1,2,\cdots. \end{equation*}
In 2000, an attempt is made by    Cvijovi\'c \&   Klinowski \cite{CK}  to use the digamma function $\psi=\Gamma'/\Gamma$ (see, e.g.,   \cite[p.136, (5.2.2)]{NIST}). They prove that
\begin{equation*}\label{(1.5)} c_0 +\pi^{-1} \psi(n+  5/4)< G_n < 1.0725+\pi^{-1} \psi(n+  5/4) ,~~n=0,1,2,\cdots ,  \end{equation*}
and
\begin{equation*}\label{(1.6)} 0.9883+\pi^{-1} \psi(n+  3 /2) < G_n < c_0+\pi^{-1} \psi(n+  3 /2) ,~~n=0,1,2,\cdots ,  \end{equation*}
where $c_0=(\gamma+4\ln 2)/\pi=1.0662\cdots$,   $\gamma=0.5772\cdots$ is the Euler constant (\cite[(5.2.3)]{NIST}).

Inequalities of this type are revisited in a  2002  paper \cite{Alzer} of Alzer. In that paper, the problem is turned   into the following: to find the largest $\alpha$ and smallest $\beta$ such that
\begin{equation*}\label{(1.7)}c_0+\pi^{-1} \psi(n+\alpha )\leq  G_n\leq  c_0+   \psi(n+\beta )~~\mbox{for~all}~ n\geq 0.  \end{equation*}
The answer is that $\alpha=5/4$ and $\beta=\psi^{-1}(\pi(1-c_0))=1.2662\cdots$, appealing to the complete monotonicity of $\Delta G_n$.

In 2009,
Zhao \cite{Zhao} starts seeking  higher terms in the bounds. A formula in \cite{Zhao},  holding for all positive integer $n$, reads
\begin{equation}\label{(1.8)} \ln (16n)+\gamma-\frac 1 {4n} +\frac 5 {192n^2} <\pi G_{n-1} <  \ln (16n)+\gamma-\frac 1 {4n} +\frac 5 {192n^2}  +\frac 3 {128n^3}. \end{equation}

Several authors have made improvements.    In a 2011 paper \cite{Mortici},
  Mortici gives an inequality of the above type involving higher order term $1/n^{5}$.  A $1/n^{7}$ term is brought in by  Granath in a recent paper \cite{Granath} in 2012.

It seems    possible to obtain sharper bounds involving terms of higher and  higher orders.  Accordingly,  difficulties  may arise.  The case by case process of taking more and more terms might  be endless.

\subsection{Asymptotic approximations}

Most  of the above inequalities can    be used to derive asymptotic approximations for $G_n$. Such approximations can also be obtained by employing integral representations, generating functions and relations with   hypergeometric functions; see, e.g.,  \cite{LLXZ}.     Indeed,
back to
 Watson \cite{Watson},  a  formula of asymptotic nature  is derived  by using a certain integral representation:
 \begin{equation}\label{(1.9)}
\begin{aligned}
     G_n& = \frac 1 \pi \ln(n+1) +\frac {\{\Gamma(n+3/2)\}^2}{\pi^2\Gamma(n+1)}\sum^{m-1}_{l=1} \frac{\{\Gamma(l+1/2)\}^2}{\Gamma(l+1)\Gamma(n+l+2)}\times \\
      &  \times \left\{ \psi(l+n+2) -\ln (n+1)+\psi(l+1)-2\psi(l+  1/2)\right\}+  O\left \{ (n+1)^{1-m}\right \}
\end{aligned}
\end{equation}for large $n$ and  positive integer $m$.
Theoretically, an asymptotic expansion can be extracted from (\ref{(1.9)}) by substituting the large-$n$ expansions of $\Gamma$ and $\psi$ into it. In fact, Waston obtains
 \begin{equation*}\label{(1.10)}
  G_n\sim \frac 1 \pi\left [   \ln(n+1)+ \gamma+4\ln 2- \frac 1 {4  (n+1)} +\frac 5 {192  (n+1)^2} +\cdots\right ], \end{equation*}
 of which (\ref{(1.8)}) is an extended version.

We skip to some very recent progress in this direction.  In
the manuscript \cite{Ismail},   Ismail,  Li and  Rahman
  derive  a complete asymptotic expansion for the Landau constants
$G_n$,   using the asymptotic sequence $n!/(n + k)!$. The approach is based on a formula of  Ramanujan, which connects the Landau constants with  a  hypergeometric function.

Several  relevant papers are worth  mentioning.   In \cite{NemesNemes}, Nemes and Nemes derive full asymptotic expansions using a formula in \cite{CK}. They also conjecture a symmetry  property of  the coefficients in the expansion.
The conjecture has been  proved by G. Nemes himself   in  \cite{Nemes}.

\noindent
\begin{prop}\label{Prop 1}
(Nemes)  Let $0 <h < 3/2$. The Landau constants $G_n$ have the following asymptotic expansions
\begin{equation}\label{(1.11)}
G_n\sim \frac 1 \pi \ln (n+h) +\frac 1 \pi (\gamma+4\ln 2 ) - \sum_{k\geq 1}\frac {g_k(h)}{(n+h)^k}\end{equation}
as $n\rightarrow +\infty$,  where the coefficients $g_k(h)$ are certain computable constants that satisfy $g_k(h)=(-1)^k g_k(3/2-h)$ for every $k\geq 1$.\end{prop}

As an important special case, Nemes \cite{Nemes} has further proved that
\begin{equation}\label{(1.12)}
\pi G_n\sim   \ln (n+3/4) + \gamma+4\ln 2  + \sum_{s=1}^\infty \frac { \beta_{2s}}{ (n+3/4)^{2s}},~~n\rightarrow \infty, \end{equation}
where the coefficients $(-1)^{s+1}\beta_{2s}$ are positive rational numbers.

The argument in  \cite{Nemes} is based on an integral representation of $G_n$ involving a  Gauss hypergeometric function in the integrand.
While in \cite{LLXZ},  the authors of the present paper study this asymptotic problem by using an entirely different approach, starting from an obvious observation that the Landau constants satisfy
a difference equation
\begin{equation}\label{(1.13)}
G_{n+1}-G_n=\left[\frac{2n+1}{2n+2}\right]^2 (G_n-G_{n-1}),~~~n=0,1,\cdots,
\end{equation}as can be seen from the explicit formula    (\ref{(1.1)}),    where  $G_{-1}:=0$.

By applying the theory of Wong and Li for
  second-order linear
difference equations \cite{WongLi1992a} to (\ref{(1.13)}),  the  general expansion  in   (\ref{(1.11)}) is obtained, and the conjecture of \cite{NemesNemes} is also confirmed.
An advantage  of this approach,
compared with the previous ones,  is that all coefficients in the expansion are given iteratively in an explicit manner.

\subsection{A question and numerical evidences}

As pointed out in \cite{LLXZ}, the case corresponding to  (\ref{(1.12)})  is numerically efficient since all odd terms in the expansion vanish.
We will find that this  expansion in terms of $n+3/4$ is even more special, both from  asymptotic  and sharper bound   points of view.

From  (\ref{(1.12)}), as suggested by the alternating signs and by numerical calculations,
there is a natural  question as follows:

\noindent {\qe\label{question 1}{Is  the error due to truncation  of  (\ref{(1.12)})  bounded in absolute value  by, and of the same sign as,  the first neglected term?
Or, more precisely,
do we have the following?
 \begin{equation}\label{(1.14)}\frac {\varepsilon_l(N)}{\beta_{2l}/N^{2l}} \in (0, 1)  ~~\mbox{for}~~n=0,1,2,\cdots~~\mbox{and}~~l=1,2,\cdots, \end{equation} where  $N=n+3/4$, and
 \begin{equation}\label{(1.15)}\varepsilon_l(N)=\pi G_n - \left \{  \ln N+\gamma+4\ln 2+\sum_{s=1}^{l-1} \frac{\beta_{2s}}{N^{2s}}\right\}.\end{equation}
}}\vskip .5cm

Recalling that  $(-1)^{s+1}\beta_{2s}$ are positive,    it is readily seen that a positive answer to  (\ref{(1.14)})  is  equivalent  to
\begin{equation}\label{(1.16)}
\varepsilon_{2k}(N)< 0~~\mbox{and}~~\varepsilon_{2k-1}(N)>0
\end{equation}for all $k=1,2,3,\cdots$ and $n=0,1,2,\cdots$.

The question reminds us of an earlier work of  Shivakumar and   Wong \cite{ShivakumarWong}, where an asymptotic expansion is obtained for the   Lebesgue  constants  associated with the polynomial interpolation at the zeros of the Chebyshev polynomials, and the error in stopping the series at any time    is shown to have the sign as, and is in absolute value less than, the first term neglected.
Similar discussion can be found in, e.g., Olver \cite[p.285]{olver1974}, on the Euler-Maclaurin formula.

\begin{figure}[h]
\begin{center}
\includegraphics[height=6cm]{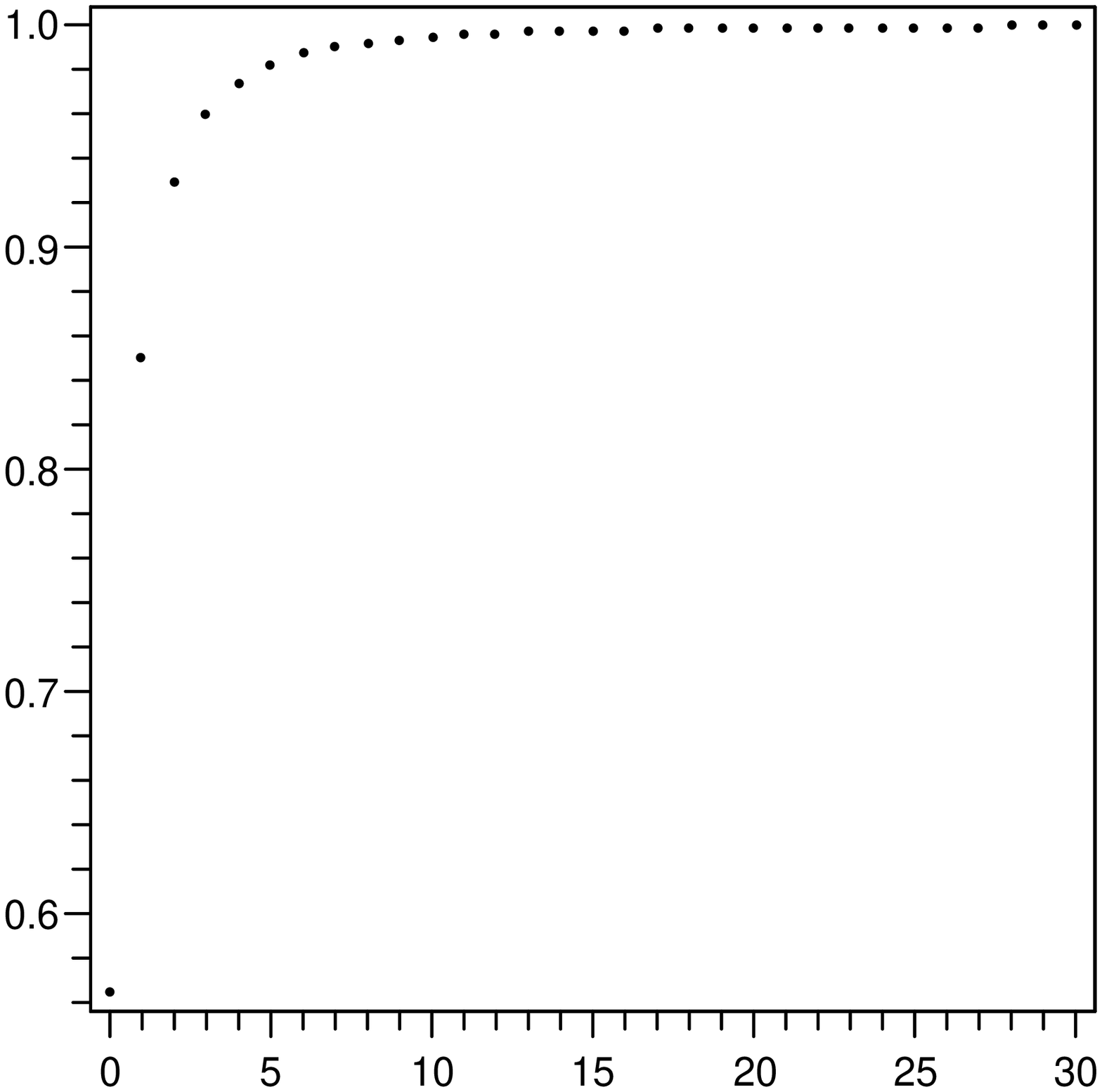}\hskip 1cm \includegraphics[height=6.12cm]{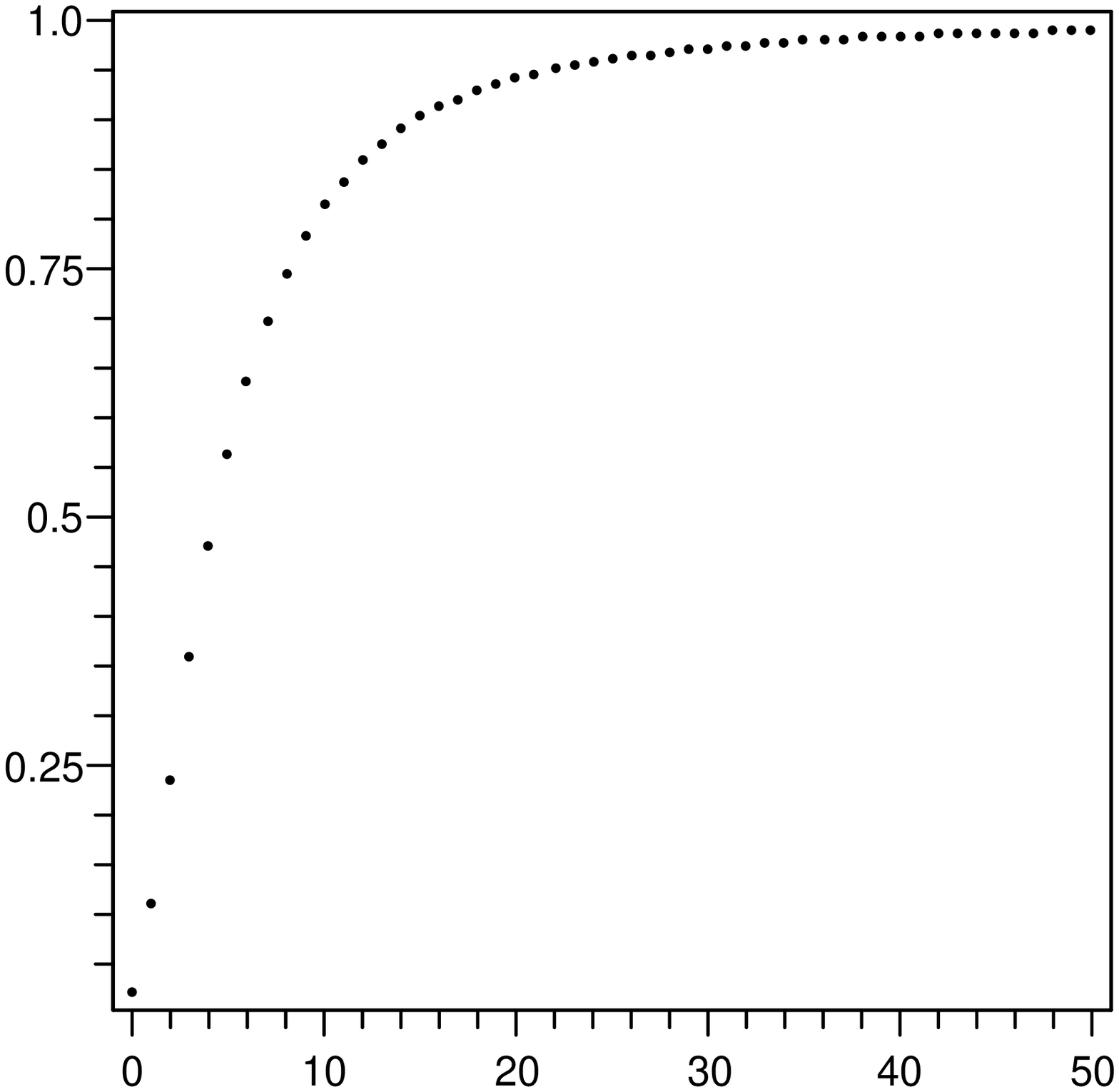}
 \caption{The function $\frac {\varepsilon_l(N)}{\beta_{2l}/N^{2l}}$.  Left: $l=2$, $n=0$-$30$.  Right: $l=16$, $n=0$-$50$.}
 \label{figure 1}
\end{center}
\end{figure}

Numerical experiments agree with (\ref{(1.14)}). The functions $\frac {\varepsilon_l(N)}{\beta_{2l}/N^{2l}}$, $N=n+\frac 3 4$, are  depicted in Figure \ref{figure 1} for the first few $n$.

\subsection{Statement of results }

In the present paper, we will justify (\ref{(1.16)}).  In fact, we will prove the following theorem.

\noindent
\begin{thm}\label{Thm 1}For $N=n+3/4$,
it holds
\begin{equation}\label{(1.17)}
(-1)^{l+1}  \varepsilon_l(N) >0\end{equation}
for  $l= 1,2,3,\cdots$ and $n=0,1,2,\cdots$,    where  $\varepsilon_l(N)$ is defined in  (\ref{(1.15)}), the coefficients $\beta_{2s}$ are determined iteratively in   (\ref{(2.4)}) below.  \end{thm}

The above theorem has direct applications both in asymptotics and sharp bounds. In asymptotic  point of view, we can obtain    error bounds which in a sense are optimal. To be precise,  we have the following.

\noindent
\begin{thm}\label{Thm 2}The error due to truncation of      (\ref{(1.12)})
is  bounded in absolute value  by, and of the same sign as,  the first neglected term for all nonnegative $n$.
That is,
\begin{equation}\label{(1.18)}
0<  (-1)^{l+1}  \varepsilon_l(N) =   \left |  \varepsilon_l(N) \right | < \frac  {\left | \beta_{2l}\right |}   {N^{2l}} = \frac {(-1)^{l+1} \beta_{2l}} { N^{2l}}\end{equation}
for  $l=0,1,2\cdots$ and $n=0,1,2,\cdots$, where $N=n+3/4$.  \end{thm}

The error bound  in (\ref{(1.18)})  is the first neglected term in the asymptotic expansion, and hence is optimal and   can not be improved.   The inequalities in (\ref{(1.18)}) can be derived from Theorem \ref{Thm 1}  by noticing that $\varepsilon_l(N)= \beta_{2l}/N^{2l}+\varepsilon_{l+1}(N)$, as can be seen from  (\ref{(1.15)}).

Another application of Theorem \ref{Thm 1} is   the construction of sharp  bounds up to arbitrary orders.

\noindent
\begin{thm}\label{Thm 3}For $N=n+3/4$,
it holds
\begin{equation}\label{(1.19)}
                   \ln N+\gamma+4\ln 2+\sum_{s=1}^{2m}\frac{\beta_{2s}}{N^{2s}}<     \pi G_n <  \ln N+\gamma+4\ln 2+\sum_{s=1}^{2k-1}\frac{\beta_{2s}}{N^{2s}}
 \end{equation}
 for all $n=0,1,2,\cdots$, $m=1,2,\cdots$, and $k=1,2,\cdots$.
 \end{thm}

The inequalities in  (\ref{(1.19)}) are understood as sharp bounds on both sides up to arbitrary orders.
In a sense, the bounds are optimal and can not be improved.

The first few
  coefficients $\beta_{2s}$ are listed  in  Table \ref{table1}, as can be evaluated via      (\ref{(2.4)}).

\noindent
\begin{table}[h]
  \centering
 \begin{tabular}{|c|c|c|c|c|c|c|c|c|c|c|}
   \hline
  $\beta_2$       & $\beta_4$              &  $\beta_6$            &  $\beta_8$                   &  $\beta_{10}                $   &$\beta_{12}$     & $\beta_{14}$     \\[0.1cm]
$\frac {11}{192}$ & $\frac{-1541}{122880}$ &$\frac{63433}{8257536}$&$\frac {-9199901}{1006632960}$&$\frac { 317959723}{17716740096}$&$\frac {-14849190321163}{281406257233920}$&$\frac {717209117969}{3298534883328}$ \\[0.1cm]

        \hline
 \end{tabular}
  \caption{The first few $\beta_{2s}$, $s=1,2,\cdots, 7$.}\label{table1}
\end{table}

\section { Proof of  Theorem \ref{Thm 1} } \indent\setcounter{section} {2}
\setcounter{equation} {0} \label{sec:2}

The proof is based on the difference equation (\ref{(1.13)}) and an approximation of the coefficients $\beta_{2s}$. To justify  Theorem \ref{Thm 1}, several  lemmas are stated, and all, except one, are proved in the present section. While the validity of  Lemma \ref{lem 2.2} is the objective of the next section.

\subsection{The coefficients $\beta_s$ in (\ref{(1.12)})}

Write the difference equation (\ref{(1.13)})   in the symmetrical form
\begin{equation}\label{(2.1)}
\left (1+\frac 1 {4N}\right )^2 w(N+1) -\left (2 +\frac 1 {8N^2}\right ) w(N)+\left (1-\frac 1 {4N}\right )^2 w(N-1)=0,
\end{equation}in which $N=n+\frac 3 4$ for $n=0,1,2,\cdots$. As mentioned earlier and as in the previous paper \cite{LLXZ}, the Landau constants  $w(N)=G_n$ solves  (\ref{(2.1)}), having an asymptotic expansion
\begin{equation}\label{(2.2)} \pi G_n\sim   \ln N+ \gamma+4\ln 2 +\sum^\infty_{s=1} \frac {\beta_s} {N^s},\end{equation}   and the coefficients  are determined by a formal substitution of (\ref{(2.2)}) into   (\ref{(2.1)}); see Wong and Li \cite{WongLi1992a}.   The following result  then follows:
\noindent
\begin{lem}\label{lem 2.1}For  $N=n+\frac 3 4$, the coefficients $\beta_s$ in   expansion (\ref{(2.2)})
 fulfill
\begin{equation}\label{(2.3)}
\beta_{2k+1}=0,~~k=0,1,2,\cdots,
\end{equation}and
\begin{equation}\label{(2.4)}
 \beta_{2k}=\frac {-1} {4k^2} \left( d_{k-1, k+1} \beta_{2k-2} + d_{k-2, k+1} \beta_{2k-4}+\cdots + d_{1, k+1} \beta_2 -d_{0, k+1}\right ),~~k=1,2,\cdots,
\end{equation}  where $d_{j, j+1}=
                      4j^2$  for $j=1,2,\cdots$,
\begin{equation}\label{(2.5)}
 d_{j, s}=                       \frac {  (2s+2j-2) \;(2s-2)! }{(2s-2j)!\; (2j-1)!} + \frac {  (2s-3)!}{8(2s-2j-2)!\;  (2j-1)!}~~\mbox{for}  ~~ s\geq j+2,
                    \end{equation}and
\begin{equation}\label{(2.6)}
 d_{0, s}= \left (\frac 1 s-\frac 1 {2s-1} \right ) +\frac 1 {16(s-1)},~~s=2,3,\cdots.
\end{equation}
In addition, it holds
\begin{equation}\label{(2.7)}
 \beta_{2k}=(-1)^{k+1}\left |   \beta_{2k}\right |,~~k=1,2,\cdots.
\end{equation}
\end{lem}\vskip .5cm

Part of this lemma ( (\ref{(2.3)}) and (\ref{(2.7)}) ) has been proved in Nemes' recent paper \cite{Nemes}. Part of it, namely (\ref{(2.3)}), and an equivalent  form of (\ref{(2.4)}),  has been proved in our earlier paper \cite{LLXZ}.  Following  Wong and Li \cite{WongLi1992a}, (\ref{(2.3)}) and (\ref{(2.4)}) can be justified by  substituting   (\ref{(2.2)}) into
(\ref{(2.1)}), expanding both sides in formal power series  of $1/N$, and equalizing the coefficients of the same powers.

It is readily seen that all $d_{j, s}>0$ for $l\geq 1$ and $s\geq l+1$, and $d_{0, s}>0$ for   $s\geq 2$.

\subsection{Analysis of  $R_l(N)$}

Here,
\begin{equation}\label{(2.8)}
R_l(N)= \left (1+\frac 1 {4N}\right )^2\varepsilon_l(N+1) -\left (2 +\frac 1 {8N^2}\right ) \varepsilon_l(N)+\left (1-\frac 1 {4N}\right )^2\varepsilon_l(N-1),
\end{equation}with  the error term  $\varepsilon_l(N)$ being  given in  (\ref{(1.15)}), and $N=n+3/4$.

There are several facts worth mentioning.  It is readily seen from \eqref{(1.15)} that $\pi G_n-(\gamma+4\ln2)$ satisfies the difference equation  \eqref{(2.1)}, and can then be removed from $R_l(N)$ in \eqref{(2.8)}. If we write $x=1/N$, then the logarithmic singularity at $x=0$ is  also  cancelled in $R_l(N)$.  Therefore,
  each $R_l(N)$ is an analytic function in $x$ for $|x|<1$. Hence the asymptotic expansion  for $R_l(N)$,  in descending powers of $N$,  is actually a convergent Taylor expansion in $x$,
\begin{equation}\label{(2.9)}
R_l(N)= \sum_{s=l+1}^\infty r_{l,s} x^{2s}, ~~ |x|<1,
\end{equation}where
\begin{equation}\label{(2.10)}
r_{l,s} = -\left ( d_{l-1, s} \beta_{2l-2} + d_{l-2, s} \beta_{2l-4}+\cdots + d_{1, s} \beta_2 -d_{0, s}\right  )
\end{equation}for $s\geq l+1$, with the leading coefficient $r_{l,l+1}=4l^2\beta_{2l}$; cf. (\ref{(2.4)}).

For later use, we   estimate     the ratio of the consecutive coefficients $\beta_{2s}$. To this end, we
introduce a sequence of positive constants
\begin{equation}\label{(2.11)}\rho_0=1, ~~\mbox{and}~\rho_l=\frac {(-1)^{l+1} \beta_{2l}} {(2l-1)!},~~l=1,2,\cdots .\end{equation}    We shall use the  following lemma and        leave its proof  to Section \ref{sec:3} below.

\noindent\begin{lem}\label{lem 2.2} It holds
\begin{equation}\label{(2.12)}
\rho_k/\rho_{k+1}\leq \frac {44} 9 \pi^2
\end{equation} for $k=0,1,2,\cdots$.
\end{lem}

Now
 we proceed to analyze $R_l(N)$ (sometimes denoted by $R_l(x)$, understanding that $x=1/N$), so as  to  show that $R_l(x)/\beta_{2l} \geq 0$ for $x\in [0, 1)$. More precisely,  we prove a much stronger result, as follows:

\noindent
\begin{lem}\label{lem 2.3} For $l=1,2,\cdots$ and $N=n+3/4$ with $n=1,2,\cdots $, we have
\begin{equation}\label{(2.13)}
(-1)^{l+1}  r_{l,s} >0, ~~ s\geq l+1 ,
\end{equation}
where $r_{l,s}$ are given   in (\ref{(2.9)}) and (\ref{(2.10)}).
\end{lem}

\noindent
{\bf{Proof:}}
The lemma can be proved by using induction with respect to $l$. Initially, we have
\begin{equation*}\label{(2.14)}
R_1(N)= \sum_{s=2}^\infty  d_{0,s}  x^{2s}.
\end{equation*}Since  $d_{0, s} >0$ for $s=2,3,\cdots$; cf. (\ref{(2.6)}), we see that  (\ref{(2.13)}) holds for $l=1$.

In view of the fact that  $\beta_2=\frac {11}{192}$; cf. Table \ref{table1}, it is readily verified that
\begin{equation*}\label{(2.15)}
R_2(N)=\sum_{s=3}^\infty r_{2,s} x^{2s} =\sum_{s=3}^\infty \left ( -\frac {11}{192} d_{1, s}+d_{0,s}\right ) x^{2s},
\end{equation*}with  all coefficients $r_{2,s}$ being negative. Indeed, in view of  (\ref{(2.5)}) and  (\ref{(2.6)}), we have
$$ r_{2,s} =-\frac {11}{192}\left (2s+\frac {2s-3} 8\right ) +\left [ \left (\frac 1 s-\frac 1 {2s-1}\right )+\frac 1 {16(s-1)}\right ] <-\frac {11}{192} \cdot (2s)+\frac 1 s <0$$
for $s\geq 3$.
 Thus (\ref{(2.13)}) holds for $l=2$.

 Similarly, we can verify (\ref{(2.13)}) for $l=3$, recalling that $\beta_4=-\frac {1541}{122880}$; cf. Table \ref{table1}.

 Assume that for $l\leq k$, it holds  $(-1)^{l+1}  r_{l,s} >0$ for  $s\geq l+1$.
 From (\ref{(2.10)}), we can write
\begin{equation}\label{(2.16)}(-1)^{k+3}  r_{k+2,s} =(-1)^{k+1}r_{k, s} +(-1)^k \left (d_{k+1, s} \beta_{2k+2}+d_{k,s}\beta_{2k}\right ) \end{equation}for $s\geq k+3 $.
 To show that  (\ref{(2.13)}) is valid for $l=k+2$, it suffices to show that
\begin{equation*}\label{(2.17)}(-1)^k \left (d_{k+1, s} \beta_{2k+2}+d_{k,s}\beta_{2k}\right )>0\end{equation*} for $s\geq k+4 $ since the validity of (\ref{(2.13)}) for $(l,s)=(k+2, k+3)$ is trivial. This is equivalent to show that
\begin{equation}\label{(2.18)}c_{k+1, s}-c_{k, s} \rho_k/\rho_{k+1} = \frac {   (-1)^k \left (d_{k+1, s} \beta_{2k+2}+d_{k,s}\beta_{2k}\right ) } {8(s-1)^2\;  (2s-3)!\;  \rho_{k+1}}  >0,\end{equation}
where $\rho_k>0$ is defined in \eqref{(2.11)}, and  $c_{k,s}=\frac {(2k-1)! d_{k,s}}{8(s-1)^2 (2s-3)!}$; cf. (\ref{(3.2)}) below.  In view of (\ref{(2.5)}), we may write
$$
c_{k,s}=\left\{\frac 1 {2(2s-2k)!}+\frac k {2(s-1) (2s-2k)!} \right\}  +\left\{ \frac 1 {64(s-1)^2(2s-2k-2)!}\right\}:=A+B.$$ For $k\geq 1$ and $s-k\geq 4$, we have
$$c_{k+1,s}\geq (2s-2k-1)(2s-2k) A+(2s-2k-3)(2s-2k-2)B\geq 56A+30B >\frac {478} 9c_{k,s}.$$
The last inequality holds since $A>8B$.

From (\ref{(2.12)}) in Lemma \ref{lem 2.2}, it is readily verified that
$$
\frac {478} 9   > \frac {44} 9 \pi^2\geq  \frac {\rho_k}{\rho_{k+1}}$$ for $k\geq 1$. Then (\ref{(2.18)}) holds for $s\geq k+4 $, and   it follows that  $(-1)^k \left (d_{k+1, s} \beta_{2k+2}+d_{k,s}\beta_{2k}\right )>0$ for $s\geq k+4 $.
Accordingly, from (\ref{(2.16)}) we see   that (\ref{(2.13)}) holds for $l=k+2$. This completes   the proof of Lemma \ref{lem 2.3}.

\subsection{Proof of Theorem \ref{Thm 1} }

Now Lemma \ref{lem 2.3} implies  that $(-1)^{l+1}  R_l(N) >0$ for all $l$ and all $N=n+3/4$.

To show that $\tilde\varepsilon_l(N):=(-1)^{l+1} \varepsilon_l (N)>0$ for all $N$, we note  first that $\tilde\varepsilon_l(N)= \frac {|\beta_{2l}|}{N^{2l}}\left\{1+O\left (\frac 1 {N^2}\right )\right \}$ as $N\rightarrow +\infty$. Hence
$\tilde\varepsilon_l (N)>0$ for $N$ large enough. Now assume that  (\ref{(1.17)}) is not true. Then there exists a finite  $M$ defined as
$$M=\max\{N= n+3/4: ~n\in \mathbb{Z}~\mbox{and}~\tilde\varepsilon_l(N)\leq 0\},$$
so that
  $M-3/4$ is a positive integer  and     $\tilde\varepsilon_l(M)\leq 0$, while $\tilde\varepsilon_l(M+1), \tilde\varepsilon_l(M+2), \cdots >0$.
For simplicity, we denote $a_N=(1+\frac 1 {4N} )^2$ and  $b_N=(1-\frac 1 {4N} )^2$.   From (\ref{(2.8)}) we have
\begin{equation*}\label{(2.19)}
a_{M+1}\tilde\varepsilon_l(M+2)=(a_{M+1}+b_{M+1}) \tilde\varepsilon_l(M+1) +  b_{M+1} \left (- \tilde\varepsilon_l(M)\right ) + (-1)^{l+1} R_l(M+1).
\end{equation*}The later terms on the right-hand side are non-negative, hence we obtain
\begin{equation*}\label{(2.20)}
a_{M+1}\tilde\varepsilon_l(M+2)\geq \left (a_{M+1}+b_{M+1}\right ) \tilde\varepsilon_l(M+1)  ,
\end{equation*}which implies that
\begin{equation}\label{(2.21)}
 \tilde\varepsilon_l(M+2)> \tilde\varepsilon_l(M+1).
\end{equation}

Using (\ref{(2.8)}) again for $N=M+2$,  we have
\begin{equation*}\label{(2.22)}
a_{M+2}\tilde\varepsilon_l(M+3)\geq (a_{M+2}+b_{M+2}) \tilde\varepsilon_l(M+2) -  b_{M+2}     \tilde\varepsilon_l(M+1) .
\end{equation*}
A combination of  the previous two inequalities gives
\begin{equation}\label{(2.23)}
\tilde\varepsilon_l(M+3)> \tilde\varepsilon_l(M+2).
\end{equation}

In general, we obtain
\begin{equation}\label{(2.24)}
\tilde\varepsilon_l(M+k+1)>  \tilde\varepsilon_l(M+k)
\end{equation}by induction. From the  equalities  in (\ref{(2.21)}),  (\ref{(2.23)}) and (\ref{(2.24)}),    we conclude that
\begin{equation}\label{(2.25)}
\tilde\varepsilon_l(M+k)>  \tilde\varepsilon_l(M+1)
\end{equation} for all $k\geq 2$. Recalling that $\tilde\varepsilon_l(N)= O(N^{-2l})$ for $N\rightarrow\infty$, letting $k\rightarrow \infty$ will give $\tilde\varepsilon_l(M+1) \leq 0$. This  contradicts the definition of $M$. Thus we have proved Theorem \ref{Thm 1}.

\section { Proof of Lemma \ref{lem 2.2}} \indent\setcounter{section} {3}\setcounter{subsection} {0}
\setcounter{equation} {0} \label{sec:3}
The idea is simple:  to approximate the coefficients     $\beta_{2s}$, and then to work out the ratio $\beta_{2s}/\beta_{2s+2}$. Yet the procedure is complicated.

 A brief outline of the proof is as follows: In Section \ref{sec:3.1}, we bring in an ordinary differential equation \eqref{(3.11)} with a specific  analytic solution $v(z)$, of which $\rho_k=\frac {(-1)^{k+1} \beta_{2k}} {(2k-1)!}$ are coefficients of the Maclaurin expansion. The function $v(z)$ is then extended, in Section \ref{sec:3.2}, and via the hypergeometric functions, to a  function analytic in the cut-strip $\{z ~|~ -4\pi<\Re z< 4\pi, ~z\not \in \{(-\infty, -2\pi]\cup [2\pi, +\infty)\}  \}$.
An integral representation is then obtained by using the Cauchy integral formula, and the integration path is deformed based on the analytic continuation procedure. In Section \ref{sec:3.3}, the integral is spilt,  approximated, and  estimated, and hence bounds for $\rho_k$  on both sides are established   in  \eqref{(3.36)}  for all $k\geq 10$. Eventually, in Section  \ref{sec:3.4}, an upper bound for $\rho_k/\rho_{k+1}$ is obtained for all non-negative integer $k$.

\subsection{Differential equation}\label{sec:3.1}
In terms of
 $\rho_s$ defined in   (\ref{(2.11)}), namely,
$\rho_0=1$ and $\rho_l=\frac {(-1)^{l+1} \beta_{2l}} {(2l-1)!}$ for $l=1,2,\cdots$,
formula (\ref{(2.4)})  can be written as
\begin{equation}\label{(3.1)}c_{l,l+1}
\rho_l -c_{l-1,l+1} \rho_{l-1} +\cdots  +  (-1)^{k}  c_{l-k, l+1} \rho_{l-k} +\cdots  +
    (-1)^{l-1}  c_{1,l+1} \rho_1 +(-1)^l    c_{0,l+1}\rho_0  =0 \end{equation}  for $l=1,2,\cdots$,
where $c_{l,l+1}= \frac 1 {2!}$, and
\begin{equation}\label{(3.2)}
c_{l-k,l+1}=\frac {(2l-2k-1)!}{(2l-1)!}\frac {d_{l-k, l+1}}{2d_{l, l+1}}= \frac {1}{2(2k+2)!} +\frac {l-k}{2l (2k+2)!}+\frac 1 {64 l^2 \cdot (2k)!}\end{equation}     for   $k=1,2,\cdots l-1$ and  $l=1,2,\cdots$. It can be verified from  (\ref{(2.6)})   that $c_{0, l+1}=\frac {d_{0, l+1}}{8l^2(2l-1)!}$ also takes the same  form, that is,   (\ref{(3.2)}) is also valid for $k=l$, $l=1,2,\cdots$.

The idea  now is to approximate $\rho_s$,  and then to estimate the ratio $\rho_{s}/\rho_{s+1}$.

Taking $l=1,2,\cdots$ in (\ref{(3.1)}), we have
  $$
  \begin{array}{l}
     a_1:=c_{1,2}\rho_1-c_{0,2}\rho_0=0, \\[0.2cm]
      a_2:=c_{2,3}\rho_2-c_{1,3}\rho_1+c_{0,3}\rho_0=0,\\[0.2cm]
      a_3:=c_{3,4}\rho_3 - c_{2,4}\rho_2+c_{1,4}\rho_1-c_{0,4}\rho_0=0,\\[0.2cm]
     \cdots\cdots.
  \end{array}$$
Summing up $\sum^\infty_{s=1} a_s x^{2s}$ gives
\begin{equation}\label{(3.3)}\rho_0\left (-c_{0,2} x^2+c_{0,3}x^4-c_{0,4} x^6+\cdots\right )+\sum^\infty_{k=1}
\rho_k x^{2k} \sum^\infty_{s=1} (-1)^{s-1} c_{k,k+s} x^{2s-2}=0.\end{equation}
In view of   (\ref{(3.2)}),   it is readily verified  by summing up the series that
\begin{equation}\label{(3.4)}-c_{0,2} x^2+c_{0,3}x^4-c_{0,4} x^6+\cdots=-\frac 1 4 +\frac {h(x)}{2} -   \int^x_0  \frac {dt_1} {t_1}\int^{t_1}_0 \frac {t h(t) dt}{16}, \end{equation}
where \begin{equation*}\label{(3.5)}h(x):=\frac {1-\cos x} {x^2}.\end{equation*}
Also we have, for $k=1,2,\cdots$,
\begin{equation}\label{(3.6)}\sum^\infty_{s=1} (-1)^{s-1} c_{k,k+s} x^{2s-2} =\frac { h(x)} 2+\frac k {x^{2k}} \int^x_0 t^{2k-1} h(t) dt - \frac 1 { x^{2k}}\int^x_0\frac {dt_1}{t_1}\int^{t_1}_0\frac{t^{2k+1}h(t)dt }{16}.\end{equation}
Substituting    (\ref{(3.4)}) and  (\ref{(3.6)}) into  (\ref{(3.3)}),   we obtain an equation
\begin{equation}\label{(3.7)}-\frac 1 4 +\frac 1 2 h(x)u(x)+\int^x_0\frac 1 2 h(t) u'(t) dt -\int^x_0\frac {dt_1}{t_1}\int^{t_1}_0\frac  t{16} h(t) u(t) dt=0,\end{equation}
where \begin{equation}\label{(3.8)}u(x):=\sum^\infty_{k=0}\rho_k x^{2k}.\end{equation}

\noindent
\begin{rem}\label{rem 1}
The existence of $u(x)$ defined above  and the validity of  (\ref{(3.3)}) can be justified by showing that $ \left |\rho_k\right |\leq M_0/\delta^{2k}$ for all positive integers $k$, $M_0$ and $\delta$ being positive constants.
Indeed, from  (\ref{(3.2)}) it is readily seen that $\left |c_{l-k, l+1} \right |\leq \frac 1 {(2k)!}$ for $k=1,2,\cdots , l$. Now we assume that $ \left |\rho_k\right |\leq M_0/\delta^{2k}$ for $k<l$, where $\delta$ is small enough such that  $2(\cosh \delta -1) <1$. Then, by using  (\ref{(3.1)}) we have
$$\frac 1 2 |\rho_l| \leq \sum^l_{k=1}  \left |c_{l-k, l+1} \right | |\rho_{l-k}| \leq \sum^l_{k=1} \frac {M_0\delta^{2k-2l}} {(2k)!} \leq \frac {M_0}{\delta^{2l}} \left (\cosh\delta -1\right ).$$
Hence we have $ \left |\rho_l \right |\leq M_0/\delta^{2l}$ by induction.
\end{rem}

Applying the operator $\displaystyle{\frac d {dx} \left ( x \frac d {dx} \right )}$ to both sides of (\ref{(3.7)}), we see that $u(x)$ solves the  second order differential equation
 \begin{equation}\label{(3.9)}\left [ x\left (\frac 1 2 h'(x) u(x)+h(x) u'(x)\right )\right ]' -\frac {xh(x) u(x)}{16} =0 \end{equation}
in a neighborhood of $x=0$, with initial conditions $u(0)=1$ and $u'(0)=0$.

In the next few steps we derive a representation of $u(x)$ for later use. First, substituting
\begin{equation}\label{(3.10)}v(x)=\sqrt{h(x)} u(x)=\frac {\sqrt 2 \sin\frac x 2} x u(x)\end{equation}
into    equation  (\ref{(3.9)})    yields
 \begin{equation}\label{(3.11)}\sin\frac x 2 v''(x)+\frac 1 2 \cos\frac x 2  v'(x) -\frac 1 {16}  \sin \frac x 2  v(x)=0  \end{equation}in a neighborhood of $x=0$, with $v(0)=1/\sqrt 2$ and $v'(0)=0$.

It is shown in Remark \ref{rem 1} that $u(x)$ is analytic at the origin. So is  $v(x)$; cf. \eqref{(3.10)}. What is more, near $x=0$, the function $v(x)$ can be represented as a hypergeometric function. Indeed,
a change of variable
\begin{equation*}\label{(3.12)}t=\frac {1-\cos x} 2=\sin^2\frac x 2\end{equation*}    turns the equation into the hypergeometric equation
\begin{equation}\label{(3.13)} t(1-t)\frac {d^2 v}{d t^2}+\left (1-\frac 3 2 t\right )\frac {dv}{dt}-\frac 1 {16} v=0.\end{equation} Taking the initial conditions into account,
it is easily verified that
\begin{equation}\label{(3.14)}v(x)=  \frac 1 {\sqrt 2} F\left (\frac 14, \frac 1 4; 1; \sin^2\frac x 2\right ) =\frac 1 {\sqrt 2} F\left (\frac 12, \frac 1 2; 1; \sin^2\frac x 4\right ) ,\end{equation} initially at $x=0$, and then analytically extended elsewhere; cf. \cite[(15.2.1)]{NIST}. The second equality follows from a   quadratic hypergeometric transformation; see \cite[(15.8.18)]{NIST}.

\subsection{Analytic continuation}\label{sec:3.2}

Well-known formulas for hypergeometric functions include
\begin{equation}\label{(3.15)}
F(a,b;c;t)=\frac {\Gamma(c)}{\Gamma(b)\Gamma(c-b)} \int^1_0  s^{b-1}(1-s)^{c-b-1}(1-st)^{-a} ds\end{equation} for $t\in \mathbb{C}\backslash [1, +\infty)$ as $\Re c>\Re b >0$; cf. \cite[(15.6.1)]{NIST}, which extends the hypergeometric function to a single-valued analytic function in  the cut-plane.

Now we proceed to consider (\ref{(3.11)}) with complex variable $z$. It is worth noting that the solution $v(z)$ we seek is an even function. So we restrict ourselves  to its analytic continuation   on the right-half plane $\Re z >0$. To this aim, we define
\begin{equation}\label{f-def} f(z) =F\left (\frac 1 2,\frac 1 2;1;\sin^2\frac z 4\right ) ~~\mbox{for}~\Re z\in (-2\pi, 2\pi)\cup (2\pi, 6\pi).\end{equation}
We see that $t=\sin^2\frac z 4$ maps the strip $-2\pi< \Re z <2\pi$ duplicately  and analytically onto the cut-plane $t\in \mathbb{C}\backslash [1, +\infty)$. The same is true for $2\pi< \Re z <6\pi$. Hence from (\ref{(3.14)})  we have the analytic continuation
$v(z)=f(z)$ for $-2\pi<\Re z<2\pi$.
  Moreover, from \eqref{f-def} it is  readily seen that the function $f(z)$, defined and analytic in  the disjoint strips, satisfies
  \begin{equation*}\label{f-periodic}f(z)=f(z-4\pi)~~\mbox{for}~~\Re z\in (2\pi, 6\pi).
\end{equation*}

Careful treatment should be brought in here, since there is a  logarithmic singularity of $v(z)$ at the boundary point  $z= 2\pi$, as we will see later, and as can be seen from the equation (\ref{(3.11)}): $z=2\pi$ is one of the nearest regular singularities, and the indicial equation has a double root $0$ there.

Next, we extend the domain of analyticity of $v(z)$ beyond the vertical line $\Re z=2\pi$.  To do so, we recall the jump along the  branch cut of the hypergeometric function
\begin{equation}\label{(3.21)}F(a,b;c;x+i0)- F(a,b;c;x-i0)  =  \frac {2\pi i \Gamma(c)}{\Gamma(a) \Gamma(b)} \frac { (x-1)^{c-a-b}} {\Gamma(c-a-b+1)} F(c-a, c-b; c-a-b+1; 1-x)
\end{equation}
for $x>1$; see  \cite[(15.2.2)-(15.2.3)]{NIST}.

Applying (\ref{(3.21)}) to $f(z)$ defined in  \eqref{f-def}, a careful calculation yields
\begin{equation*}\label{f-jump-upper}f(2\pi -0+iy)-f(2\pi +0+iy)=2i f(iy )~~\mbox{for}~~y=\Im z >0,
\end{equation*}
where use has been made of the fact that
  $t=\sin^2 \frac {2\pi-0+iy} 4$, $y=\Im z>0$ corresponds to the upper edge of $[1, \infty)$, while $t=\sin^2 \frac {2\pi+0+iy} 4=\sin^2 \frac {-2\pi+0+iy} 4$, $y=\Im z>0$  corresponds to the lower edge, and that $t=\sin^2 \frac z 4$ maps   the upper imaginary axis to  $(-\infty, 0)$.
Similarly we have
\begin{equation*}\label{f-jump-lower}f(2\pi -0+iy)-f(2\pi +0+iy)=-2i f(iy )~~\mbox{for}~~y=\Im z <0.
\end{equation*}

Summarizing the above, we have an analytic function in the cut-strip $\{z ~|~ 0\leq \Re z< 4\pi, ~z\not \in  [2\pi, +\infty)  \}$, defined as follows:
\begin{equation}\label{(3.25)}
v(z)=\left\{
\begin{array}{ll}
  f(z), & 0\leq \Re z<2\pi, \\[.1cm]
  f(z)  +2if(z-2\pi), & 2\pi<\Re z <4\pi,~ \Im z >0, \\[.1cm]
f(z)  -2if(z-2\pi), & 2\pi<\Re z <4\pi,~ \Im z <0.
\end{array}\right .
\end{equation}The value at $\Re z=2\pi$, $\Im z\not=0$ is obtained by taking limit.

At last, we confirm that there is a logarithmic singularity at $z=2\pi$,   a regular singularity of (\ref{(3.11)}). Indeed,  following the derivation from (\ref{(3.11)}) to (\ref{(3.13)}), we see   that all solutions  to (\ref{(3.11)}) takes the form
$$A w_1(z) +B \left \{ w_1(z)\ln \left (\sin^2\frac z 2\right ) + w_2(z)\right\}=A w_1(z) +B \left \{2 w_1(z)\ln \left (z -2\pi \right ) + \tilde w_2(z)\right\};$$ cf. Wong \cite[(2.1.24)]{Wong2010},
 where $w_2(z)$ and  $\tilde w_2(z)$ are  specific single-valued analytic functions at $z=2\pi$,  $A$ and $B$ are constants, and $w_1(z)= f(z-2\pi)$ is an analytic solution of (\ref{(3.11)}) at $z=2\pi$. The function $v(z)$ in  (\ref{(3.25)}) is such a solution, and, what is more, with $B=-\frac 1 \pi$, as can be seen by comparing the jumps along $(2\pi, 3\pi)$. Accordingly, we may write
 \begin{equation}\label{(3.26)}
v(z)=v_A(z) -\frac 2 \pi f(z-2\pi)\ln \left (2\pi -z  \right )
\end{equation} for $0< \Re z <4\pi$, with  $v_A(z)$ being  analytic in  the strip, and the branch of  logarithm being  chosen as $\arg (2\pi -z)\in (-\pi, \pi)$.

\subsection{Approximation of $\rho_s$}\label{sec:3.3}

Now we are in a position to approximate $\rho_s$.

To begin, we mention several known facts.

From  (\ref{(3.25)}), and that $v(z)$ is even, we know that the function $v(z)$ is now extended analytically to the cut-strip $$\{z ~|~ -4\pi<\Re z< 4\pi, ~z\not \in \{(-\infty, -2\pi]\cup [2\pi, +\infty)\}  \}.$$

Again, from \eqref{f-def}, (\ref{(3.25)})  and (\ref{(3.26)}), and that the hypergeometric function $F(a,a;c;t)$ behaviors of the form $O(|t|^{-a}\ln |t|)$ at infinity; see  \cite[(15.8.8)]{NIST},
 we conclude that
$$ v(z)=O\left ( e^{-\frac 1 4 |\Im z|}\ln |\Im z |\right )$$ as $\Im z \rightarrow \infty$,  holding  uniformly in $|\Re z |\leq 4\pi -\delta$ for arbitrary positive $\delta$.
Thus  $v(z)$ displays an exponentially small decay at infinity in each sub-strip.

\begin{figure}[h]
\begin{center}
\includegraphics[height=6cm]{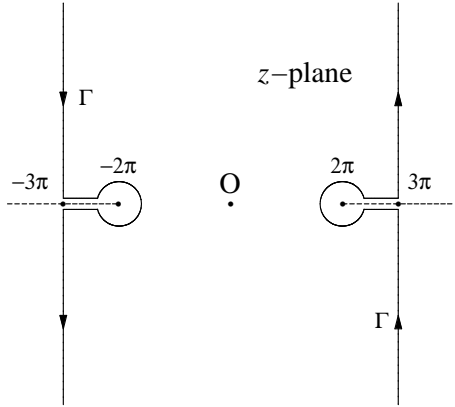}
 \caption{The deformed contour $\Gamma$: the oriented curve.}
 \label{figure 2}
\end{center}
\end{figure}

Hence, we can derive from (\ref{(3.8)}) and (\ref{(3.10)}) that
\begin{equation}\label{(3.27)}\rho_k=\frac 1 {2\pi i}  \oint z^{-2k-1} u(z)dz =\frac 1 {\sqrt 2\pi i}  \oint \frac {\frac z 2}{\sin \frac z 2}\frac { v(z)  dz} { z^{2k+1} }
  =\frac 1 {\sqrt 2\pi i} \int_\Gamma  \frac {\frac z 2}{\sin \frac z 2}\frac { v(z)  dz} { z^{2k+1} },\end{equation}
where initially the integration path  is a small circle encircling the origin, and then being   deformed to the  contour  $\Gamma$   illustrated in Figure \ref{figure 2}.

From the symmetry property of $v(z)$, we need only evaluate and estimate the integral on the right-half of $\Gamma$. We split the integral in (\ref{(3.27)}) into three integrals, namely,
\begin{equation}\label{(3.28)}\frac {\pi i} {\sqrt 2}\rho_k=   \int_{\Gamma_v}  \frac {\frac z 2}{\sin \frac z 2}\frac { v(z)  dz} { z^{2k+1} }
 +  \int_{\Gamma_l}  \frac {\frac z 2}{\sin \frac z 2}\frac { v_A(z)  dz} { z^{2k+1} }-  \int_{\Gamma_l}  \frac {\frac z \pi f(z-2\pi)\ln (2\pi-z)  dz}{\sin \frac z 2 ~ z^{2k+1} }:=I_v+I_a+I_l     ,
\end{equation} where $\Gamma_v$ is the right-half vertical part $\Re z=3\pi$, and $\Gamma_l$ is the remaining  right-half $\Gamma$, consisting of an circular part around $z=2\pi$, and a pair of horizontal line segments joining the circle with the vertical line; compare Figure \ref{figure 2} for the curves and the orientation.

We will show that the main contribution to $\rho_k$, when $k$ is not small, comes from $I_l$. Estimates will be obtained with full details. We do the calculation case by case.

\vskip .5cm
\noindent
{\bf{Estimating $I_v$}}:

First, we estimate $v(z)$   for $\Re z=3\pi$.   It is readily seen that     $t=\sin^2\frac {3\pi+iy} 4=\frac {1-i\sinh y} 2$ and $\sin^2\frac {\pi+iy} 4=\frac {1+i\sinh y} 2$, where $y=\Im z$. Hence $$\arg \left\{(1-st)^{-\frac 1 2}\right \}\in (-  \pi/ 2, 0)~\mbox{for}~y\in (0, \infty),~~\mbox{and}~\left | (1-st)^{-\frac 1 2}\right |\leq \left (1-\frac s 2\right )^{-\frac 1 2}~\mbox{for}~0\leq s \leq 1.$$  Accordingly, from (\ref{(3.15)}) and  \eqref{f-def} we have
$$ \left | f(3\pi+iy)\right | \leq \frac 1 \pi \int^1_0 s^{-\frac 1 2} (1-s)^{-\frac 1 2} \left (1-s/2\right )^{-\frac 1 2} ds:=C_{f,v}=1.1803\cdots ~~\mbox{for}~~y> 0. $$ Noting that $\arg f(3\pi +iy) \in (-\pi/2, 0)$ for $y>0$, and that $f(\pi +iy)=\overline{f(3\pi +iy)}$,
Substituting  all above into    (\ref{(3.25)}), we have
$$|v( 3\pi +iy)| \leq   \sqrt 5  \left | f(3\pi+iy)\right |\leq     M_v:=\sqrt 5\;  C_{f,v} =2.6393\cdots ,~~y>0.$$ Similar argument shows that the equality holds for $y\in \mathbb{R}$.

Hence, we conclude from    (\ref{(3.25)})
 and (\ref{(3.28)}) that
\begin{equation}\label{(3.29)}
|I_v |\leq \frac {1} 2 \int^\infty_{-\infty} \frac {|v(  3\pi+iy)|\; dy}{\cosh\frac y 2 ~ | 3\pi+iy|^{2k}}\leq \frac {M_v} {2} \frac 1 {(3\pi)^{2k}}\int^\infty_{-\infty}  \frac {dy} {\cosh\frac y 2}
=\frac { \pi M_v} {(3\pi)^{2k}}=\frac { 8.2916\cdots } {(3\pi)^{2k}}.
\end{equation}

\vskip .5cm
\noindent
{\bf{Evaluating  $I_a$}}:

This time, $v_A(z)$ is analytic in a domain containing $\Gamma_l$. It is readily seen that only the simple pole $z=2\pi$ contributes to $I_a$. More precisely, applying the residue theorem we have
$$ I_a  =\frac {2\pi i v_A(2\pi)} {(2\pi)^{2k}}.$$
Here
$v_A(2\pi)$ can be obtained by substituting $z=2\pi-\varepsilon$ into  (\ref{(3.25)}) and (\ref{(3.26)}), and   taking limit
$$ v_A(2\pi)=\lim_{\varepsilon\rightarrow 0^+}\left\{ \frac 1 {\pi} \int^1_0 s^{-\frac 1 2} (1-s)^{-\frac 1 2} (1-st)^{-\frac 1 2} ds  +  \frac 2 \pi f(-\varepsilon) \ln\varepsilon\right \},$$where $t=\sin^2\frac z 4=\cos^2\frac   \varepsilon 4$; cf. (\ref{(3.15)}) and \eqref{f-def}.

A careful calculation leads to  \begin{equation}\label{(3.30)}
v_A(2\pi)  =\frac 1 \pi\int^1_0 \frac {s^{-\frac 12 }-1}{1-s } ds  +       \frac 1 \pi\int^1_0 \frac {s^{-\frac 1 2 }-1}{1-s } ds+\frac 4 \pi \ln 2=  \frac 8 \pi \ln 2    .
\end{equation}
Indeed,  (\ref{(3.30)}) can be derived as follows:
$$\frac 1 {\pi} \int^1_0\left ( s^{-\frac 1 2} -1\right )(1-s)^{-\frac 1 2} (1-st)^{-\frac 1 2} ds \longrightarrow \frac 1 \pi\int^1_0 \frac {s^{-\frac 12 }-1}{1-s } ds ~~\mbox{as}~~t\rightarrow 1^-$$ by applying the dominated convergence theorem. Also,
$$\frac 1 {\pi} \int^1_0   (1-s)^{-\frac 1 2} (1-st)^{-\frac 1 2} ds = \frac 1 { \pi} \int^1_0 \frac {\tau ^{-\frac 1 2 } }{1-t\tau  } d\tau $$ by making change of variable $\tau=\frac {1-s}{1-ts}$. Finally,
$$\frac 1 {\pi} \int^1_0 \left\{s ^{-\frac 1 2 }-1 \right \} \frac 1 {1-ts  } ds  \longrightarrow  \frac 1 { \pi} \int^1_0 \frac {s^{-\frac 1 2 }-1}{1-s } ds~~\mbox{as}~~t\rightarrow 1^- $$   by applying the dominated convergence theorem one more time. For $t<1$, it further  holds   $$\frac 1 {\pi} \int^1_0 \frac { 1 }{1-ts  } ds=-\frac 1{\pi t} \ln (1-t).$$
Adding up all these gives   (\ref{(3.30)}).
Thus we have
\begin{equation}\label{(3.31)}
 I_a = \frac {16\ln 2} {(2\pi)^{2k}}i = \frac {11.0903\cdots} {(2\pi)^{2k}}i .
\end{equation}

\vskip .5cm
\noindent
{\bf{Evaluating  $I_l$}}:

Now we turn to the last integral $I_l$ in (\ref{(3.28)}). First, we observe that $\displaystyle{\frac {f(z-2\pi)}{\sin\frac {z-2\pi} 2}- \frac {f(0)}{ \frac {z-2\pi} 2} } $ is analytic in the strip $\Re z\in (0, 4\pi)$. Hence the circular part of $\Gamma_l$ collapses in the following integral to give
$$\frac 1 \pi \int_{\Gamma_l} \left\{ \frac {f(z-2\pi)}{\sin\frac {z-2\pi} 2}- \frac {f(0)}{ \frac {z-2\pi} 2}\right \} \frac {\ln (2\pi -z)}{z^{2k}} dz= -2i \int_{2\pi}^{3\pi}  \left\{ \frac {f(z-2\pi)}{\sin\frac {z-2\pi} 2}- \frac {f(0)}{ \frac {z-2\pi} 2}\right \} \frac {dz}{z^{2k}}  . $$
Accordingly we have
$$\left | \frac 1 \pi \int_{\Gamma_l} \left\{ \frac {f(z-2\pi)}{\sin\frac {z-2\pi} 2}- \frac {f(0)}{ \frac {z-2\pi} 2}\right \} \frac {\ln (2\pi -z)}{z^{2k}} dz\right |\leq \frac {4\pi M_f}{2k-1} \frac 1 {(2\pi)^{2k}},$$
where $$ M_f\geq  \max_{2\pi \leq z \leq 3\pi} \left |\frac {f(z)}{\sin\frac {z-2\pi} 2}- \frac {f(2\pi)}{ \frac {z-2\pi} 2}\right |.$$
We give a rough estimate for this maximum value. Recalling  that $f(z)=F\left (\frac 1 2, \frac 12; 1; \sin^2\frac {z} 4\right )$ and $f(0)=1$, noting that  the Maclaurin  expansion of  $F\left (\frac 1 2, \frac 12; 1; t \right )$ has all positive coefficients (cf. \cite[(15.2.1)]{NIST}),  and that $\frac {d}{dt} \left ( \frac {F\left (\frac 1 2, \frac 12; 1; t \right )-1}{\sqrt t}\right )>0$        for $t\in (0, 1)$ via term-by-term differentiation, we see that $\frac {f(z-2\pi)-1 }{\sin\frac {z-2\pi} 2}=\frac 1 {2\cos\frac {z-2\pi} 4} \frac {F\left (\frac 1 2, \frac 12; 1; t \right )-1}{\sqrt t}$ is monotone increasing for
 $ 2\pi<   z  \leq 3 \pi$, or correspondingly,  $t= \sin^2\frac {z-2\pi} 4\in (0, 1/2]$. Therefore, we   have
$$ 0< \frac {f(z-2\pi)-1 }{\sin\frac {z-2\pi} 2} \leq   F \left (\frac 1 2, \frac 12; 1; \frac 1 2 \right ) -1  = \frac {\sqrt\pi} {\Gamma(3/4)\Gamma(3/4)} -1; $$
see \cite[(15.4.28)]{NIST}. Also, one has
$$0< \frac 1 {\sin x}-\frac 1 x =\frac {x-\sin x} {x\sin x}\leq \frac {x^3}{6x\sin x} \leq \frac {\pi^2} {24}~~\mbox{for}~~0<x\leq \frac \pi 2.$$
Hence an appropriate   choice of $M_f$ is
$$M_f=\frac {\sqrt\pi} {\Gamma(3/4)\Gamma(3/4)} -1 + \frac {\pi^2} {24}=0.5915\cdots.$$

We proceed to evaluate the crucial part
$$\frac 2 \pi \int_{\Gamma_l}   \frac {\ln (2\pi-z)}{z-2\pi} \frac {dz}{z^{2k}},$$
of which the integrand has a pole $z=2\pi$, coinciding  with the logarithmic singularity.
To treat such kind of singularities, we appeal to the idea of
 Wong and Wyman \cite{WongWyman}.

 For simplicity we re-scale the variable $\tau=\frac {z-2\pi} {2\pi}$, which turns the integral into
 $$\frac 2 {\pi (2\pi)^{2k}} \int_{\frac 1 2}^{(0-)}   \frac {\ln (2\pi)+\ln(-\tau)}{\tau } \frac { d\tau} { e^{2k\ln (1+\tau)}}
 =-\frac {4i\ln 2\pi}{(2\pi)^{2k}} + \frac 2 {\pi(2\pi)^{2k}}   \int_{\frac 1 2}^{(0-)}   \frac { \ln(-\tau)}{\tau } \frac {d\tau }{e^{2k\ln (1+\tau)}}
  ,$$
 where the integration path  starts and ends at $\tau=\frac 1 2$, and encircles the origin clockwise: the loop is a re-scaled  version of $\Gamma_l$.   The branch of $\ln (1+\tau) $ is chosen  so that $\ln (1+\tau)>0 $ for $\tau >0$.

 To extract the main contribution, we further split the exponent  in the integrand.  Indeed, we see that
$$ \left |    \frac 2 {\pi(2\pi)^{2k}}   \int_{\frac 1 2}^{(0-)}   \frac {e^{2k(\tau -\ln (1+\tau)) }
-1}{\tau } e^{-2k\tau} \ln(-\tau)d\tau \right | \leq \frac 2 { (2\pi)^{2k}} \max_{0\leq \tau <1/2} \left |\varphi(\tau )\right |\leq \frac {2M_\varphi} { (2\pi)^{2k}},$$
 where the   function $$\varphi(\tau )=\frac {e^{2k(\tau -\ln (1+\tau)) }
-1}{\tau } e^{-2k\tau}$$ is analytic in a neighborhood of $[0, 1/2]$,  thus the integration path collapses to the lower and upper edges of $[0, 1/2]$, and its bound $M_\varphi$   can be obtained by noticing that
$$0\leq \tau -\ln(1+\tau )\leq \frac {\tau^2} 2~~\mbox{for}~~\tau >0,~~~\mbox{and}~~-\ln(1+\tau)\leq -\frac 2 3\tau ~~\mbox{for}~~ 0\leq \tau \leq \frac 1 2.$$

Now for $0\leq \tau \leq \frac 1 {\sqrt k}$, one has
$$0\leq \varphi(\tau)\leq  \frac {e^{k \tau^2}
-1}{\tau } e^{-2k\tau}    \leq (e-1) (k\tau)e^{-2(k\tau)} \leq \frac {e-1}{2e}. $$
 While for $ \frac 1 {\sqrt k}\leq \tau \leq \frac 1 2$,
 $$0\leq \varphi(\tau)\leq  \frac {e^{-2k \ln (1+\tau) }}\tau \leq \frac {e^{-\frac 4 3 k\tau}}\tau \leq \left . \frac {e^{-\frac 4 3 k\tau}}\tau \right |_{\tau=\frac 1 {\sqrt k}}=\sqrt k e^{-\frac  4 3\sqrt k} \leq \frac 3 {4e}
 . $$ Thus we may chose
 $$M_\varphi=\max\left ( \frac {e-1}{2e},    \frac 3 {4e}\right ) =\frac {e-1}{2e},$$ which does not depend on $k$.

 The remaining piece would turn out to be of the most significance. Indeed, a change of variable $s=2k\tau$ makes
\begin{equation*}\label{(3.32)}  \begin{array}{rl}
      \displaystyle{   \frac 2 {\pi(2\pi)^{2k}}   \int_{\frac 1 2}^{(0-)}   \frac {\ln (-\tau)  }
 \tau   e^{-2k\tau}  d\tau }&=  \displaystyle{    \frac 2 {\pi(2\pi)^{2k}}   \int_{k}^{(0-)} \frac {\ln(-s)-\ln 2k} s e^{-s} ds }  \\[.4cm]
        &  =\displaystyle{ \frac {4i\ln 2k} {(2\pi)^{2k}}+  \frac 2 {\pi(2\pi)^{2k}}   \int_{k}^{(0-)} \frac {\ln(-s)} s e^{-s} ds .}
     \end{array}
 \end{equation*}

   The last integral can be approximated  by   extending  the integration path to $\displaystyle{\int_{+\infty}^{(0-)}}$, with an error bounded by  $\displaystyle{\frac {4 e^{-k}}   {k(2\pi)^{2k}}}$.
Recalling the integral representation
$$\frac 1 {\Gamma(z)}=\frac 1 {2\pi i} \int^{(0-)}_{+\infty} (-s)^{-z} e^{-s}ds$$ for all finite $z$; cf. \cite[(5.9.2)]{NIST} or  Wong and Wyman \cite{WongWyman}, taking derivative with respect to  $z$ and setting $z=1$, we obtain
$$   \int^{(0-)}_{+\infty}  \frac {\ln(-s)} s e^{-s} ds= \frac {2\pi i\Gamma'(1)}{\Gamma^2(1)} =-2\pi i \gamma,$$
 where $\gamma=0.5772\cdots$ is   Euler's constant; see \cite[(5.4.12)]{NIST}.

Hence we can write
\begin{equation}\label{(3.33)} I_l=  \frac {4i\ln 2k} {(2\pi)^{2k}}+  \frac {(-4\gamma-4\ln2\pi)i+\delta_{l,k}} { (2\pi)^{2k}},
 \end{equation}
with
$$|\delta_{l,k}|\leq \frac {4\pi M_f}{2k-1} + 2M_\varphi + \frac {4e^{-k}} k .$$

\vskip .3cm

Now we   sum up the calculations and estimations  made above.  Substituting (\ref{(3.29)}), (\ref{(3.31)}) and  (\ref{(3.33)}) into (\ref{(3.28)}), we have the approximation
\begin{equation}\label{(3.34)}
\frac \pi {\sqrt 2}
\rho_k= \frac {4 \ln 2k} {(2\pi)^{2k}}+  \frac {\left (11.0903\cdots -4\gamma-4\ln2\pi\right )+\delta_{k}} { (2\pi)^{2k}}, \end{equation} with
 \begin{equation*}\label{(3.35)}
|\delta_{k}| \leq   \pi M_v \left ( \frac 2 3\right )^{2k}  +  \frac {4\pi M_f}{2k-1}  +2M_\varphi + \frac {4e^{-k}} k  < 1.0259  \end{equation*} for $k\geq 10$. Hence  we obtain from (\ref{(3.34)}) that
\begin{equation}\label{(3.36)}
\frac {4 \ln 2k  + 0.4041\cdots  } {(2\pi)^{2k}} \leq
\frac \pi {\sqrt 2}
\rho_k\leq  \frac {4 \ln 2k  +2.4559\cdots  } {(2\pi)^{2k}}  \end{equation}for $k\geq 10$.

\subsection{The ratio}\label{sec:3.4}

For $k\geq 10$, it is readily verified that
$$ \frac 8 9\ln(2k+2) \geq 2.7475\cdots >2.4559-\frac {11} 9 \times 0.4041\cdots.$$ Therefore, it follows  from  (\ref{(3.36)})     that

\begin{equation*}\label{(3.37)} \rho_k/\rho_{k+1}\leq \frac {44} 9 \pi^2 \end{equation*}for $k\geq 10$.

It is easily verified that the inequality holds for all $k\geq 0$ by numerical evaluation of the first few $\rho_k$ for $k=0,1,\cdots 9$, as can be seen from  Table \ref{table2}.

\begin{table}[h]
  \centering
 \begin{tabular}{|c|c|c|c|c|c|c|c|c|c|c|}
   \hline
  $k$                & $0$        &  $1$       &    $2$    &    $3$    &$4$        & $5$     & $6$    & $7$     &    $8$  &  $9$    \\
$\rho_k/\rho_{k+1}$  & 17.46      & 27.41      &32.65      &35.30      &36.67      &37.41    & 37.86  & 38.15   &  38.36  &  38.51  \\[0.1cm]
$\frac {44} 9\pi^2$  & 48.25      & 48.25      &48.25      &48.25      &48.25      &48.25    &48.25   &48.25    &48.25    &48.25    \\
        \hline
 \end{tabular}

  \caption{The first few $\rho_k$, $k=0,2,\cdots, 9$.}\label{table2}
\end{table}

A by-product is that  $\rho_k/\rho_{k+1} \rightarrow 4\pi^2$ as $k\rightarrow \infty$. Another byproduct of (\ref{(3.36)})  and  Table \ref{table2} is that      $\rho_k>0$ for all $k$, and thus $\beta_{2k}$ having  alternative signs, as stated in   (\ref{(2.7)}).

\section { Discussion } \indent\setcounter{section} {4}
\setcounter{equation} {0} \label{sec:4}

We discuss very briefly a   conjecture of     H. Granath, of which we were not aware   until  we almost  finish writing the present paper. In   \cite{Granath},  Granath derives an asymptotic expansion
\begin{equation}\label{(4.1)}
\pi G_{n-1}\sim  \ln(16n)+\gamma +\sum_{k=1} ^\infty  \frac{a_k}{(16n)^k}, ~~n\rightarrow\infty,   \end{equation}  where $a_k$ are `effectively computable' constants   but not explicitly given, except for the first few.
The author shows interest in seeking  sharp  bounds of arbitrary orders. Indeed, denoting
\begin{equation}\label{(4.2)}
A_m(n)=\ln(16n)+\gamma +\sum_{k=1} ^m  \frac{a_k}{(16n)^k} ,  \end{equation}
Granath proves that
$A_5(n) <\pi G_{n-1} < A_7(n)$ and states that   $A_9(n) <\pi G_{n-1} < A_{11}(n)$,   for all positive $n$. It is  then  conjectured that (the   sign in \cite[(10)]{Granath} is apparently wrong)
\begin{equation}\label{(4.3)}
(-1)^{\frac {m(m+1)} 2} \left (  \pi G_{n-1}-A_m(n)\right ) <0   \end{equation}for all  $n=1,2,\cdots$ and $m=0,1,2,\cdots$.

The existence of    (\ref{(4.1)}) has also  been justified in \cite{LLXZ} and \cite{NemesNemes}. It is easily seen  that in terms of $\beta_{s}$ in   (\ref{(2.3)}) and (\ref{(2.4)}),  the coefficients can be written as
\begin{equation}\label{(4.4)}
a_k=4^k\left [  -\frac 1 k +\sum^k_{s=1} \frac { (k-1)! 4^s \beta_s  }{(s-1)!(k-s)!} \right ], ~~k=1,2,\cdots;  \end{equation}
see also \cite{LLXZ} for an   iterative formula.  So far, numerical experiments  agree with  (\ref{(4.3)}). A proof of it   might be found  either by following the steps, or by using the results, of the present paper.

Then, a natural question  arises:

\noindent {\qe\label{question 2}{Considering the general expansion in   (\ref{(1.11)}), for what $h$ we have the ``best'' approximation in the sense of  Theorem \ref{Thm 1} and (\ref{(4.3)})?
}}

The case $h=3/4$ is what we have been dealt with in the present paper.  Very likely the expansion  (\ref{(1.11)}) for  $h=1/2$ and $h=1$ would  turn out to be the ``best''.

As mentioned earlier, the coefficients    $\beta_{2s}$  of the expansion    (\ref{(1.12)}) can be obtained iteratively via  (\ref{(2.4)}). We complete the paper by giving a couple of
 alternative representations    for   $\beta_{2s}$.  For example,
 it can be shown  that
 \begin{equation}\label{(4.5)}
 \beta_{2l} =   \frac { (-1)^{l+1}}{ 2^{2l} (l!)^2 }     \left |
                    \begin{array}{cccccc}
                     d_{0,2}  &  d_{0,3}  &  \cdots & d_{0,l-1} & d_{0, l} & d_{0, l+1} \\
                     d_{1,2}  &  d_{1,3}  &  \cdots & d_{1,l-1} & d_{1, l} & d_{1, l+1} \\
                     0        &  d_{2,3}  &  \cdots & d_{2,l-1} & d_{2, l} & d_{2, l+1} \\
                     \vdots   &  \vdots   &  \ddots & \vdots    & \vdots   & \vdots   \\
                      0       & 0         & \cdots  &d_{l-2,l-1}&d_{l-2,l} &d_{l-2,l+1} \\
                      0       & 0         & \cdots  &   0       &d_{l-1,l} &d_{l-1,l+1} \\
                    \end{array}
                  \right |, ~~~l=1,2,\cdots.\end{equation}Indeed, taking $k=1,2,\cdots, l$  in (\ref{(2.4)}), we have
a linear algebraic system with unknowns $\beta_2$, $\beta_4$, $\cdots$, $\beta_{2l}$.  Solving the system gives  (\ref{(4.5)}).

Also, a combination of (\ref{(2.11)}),    (\ref{(3.14)}) and     (\ref{(3.27)}) yields  an   integral representation
\begin{equation}\label{(4.6)}
\beta_{2l}=\frac {(-1)^{l+1} (2l-1)!}  {4\pi i} \oint \frac z {\sin\frac z 2} \frac {F\left (   1/ 4 ,  1/ 4 ; 1; \sin^2\frac z 2\right ) dz}{z^{2l+1}}  , ~~l=1,2,\cdots,\end{equation}
where the integration path is a small circle encircling  the origin counterclockwise. Of course, (\ref{(3.34)}) can be interpreted  as
\begin{equation}\label{(4.7)}
\beta_{2l}=\frac {(-1)^{l+1} (2l-1)!}  {\pi (2\pi)^{2l}}\left\{4\sqrt 2\ln(2l) +O(1) \right\}  \end{equation} for large $l$.  Results can be obtained from such approximations. For example,  the expansion (\ref{(1.12)}) is divergent; compare \cite[Theorem 3]{Granath}.

\section*{Acknowledgements}
The authors are grateful to Prof. R. Wong for bringing the problem into their attention.
The authors thank the anonymous referees for their carefully reading of the manuscript and for the  valuable suggestions and comments. One referee suggested the using of a quadratic    transformation of the hypergeometric functions which makes the proof of Lemma \ref{lem 2.2}    more  natural and simplified. The other referee provided many constructive suggestions and corrections which have much  improved the readability of the  manuscript.

The work of Yutian Li was supported in part by the HKBU Strategic Development Fund,
a start-up grant from Hong Kong Baptist University,
and a grant from the Research Grants Council of the Hong Kong Special Administrative Region, China (Project No. HKBU 201513).
The work of Saiyu Liu was supported in part by   Hunan  Natural Science Foundation under grant number 14JJ6030, and by the National Natural Science Foundation of China under grant number  11326082.
The work of Shuaixia Xu  was supported in part by the National
Natural Science Foundation of China under grant number
11201493, GuangDong Natural Science Foundation under grant number S2012040007824, and the Fundamental Research Funds for the Central Universities under grand number 13lgpy41.
 Yuqiu Zhao  was supported in part by the National
Natural Science Foundation of China under grant  numbers 10471154 and
10871212.

\end{document}